\documentclass[10pt]{article}
\usepackage{amsmath}
\usepackage{amssymb}
\usepackage{amscd}
\usepackage{latexsym}
\usepackage{theorem}
\usepackage{doublespace}
\usepackage{epsfig}
\usepackage{psfrag}

\newtheorem{theorem}{Theorem}[section]

\newtheorem{proposition}[theorem]{Proposition}
\newtheorem{corollary}[theorem]{Corollary}

\newtheorem{conjecture}[theorem]{Conjecture}

{\theorembodyfont{\rmfamily}
\theoremstyle{plain}
\newtheorem{definition}[theorem]{Definition}
\newtheorem{example}[theorem]{Example}

\newtheorem{remark}[theorem]{Remark}
}

\renewcommand{\Re}{\operatorname{Re}}

\newcommand{\C}{{\mathbb{C}}}
\newcommand{\Z}{{\mathbb{Z}}}

\newcommand{\R}{{\mathbb{R}}}
\newcommand{\Zt}{{\mathbb{Z}_2}}

\newcommand{\la}{\lambda}
\renewcommand{\a}{\alpha}

\newcommand{\F}{\mathcal{F}}
\newcommand{\Sala}{\rm{Sal}(\A)}

\newcommand{\subs}{\subseteq}

\newcommand{\hookto}{{\hookrightarrow}}

\newcommand{\<}{\left<}
\renewcommand{\>}{\right>}

\newcommand{\hs}{\hspace{3pt}}

\newcommand{\Sp}{S^+}
\newcommand{\Sm}{S^-}
\newcommand{\A}{\mathcal{A}}
\newcommand{\conea}{\mathcal{C}\mathcal{A}}
\newcommand{\CA}{\mathcal{M}_{\R}(\A)}
\newcommand{\MA}{\mathcal{M}(\A)}

\newcommand{\hz}{H^*_{\Zt}}

\newcommand{\Hip}{H_i^+}
\newcommand{\Him}{H_i^-}

\newcommand{\vgzt}{VG(\A;\Zt)}

\newcommand{\qed}{\hfill \mbox{$\Box$}\medskip\newline}
\newenvironment{proof}{\noindent {\bf Proof:}}{\qed \par}

\begin{document}
\begin{spacing}{1.1}

\noindent
{\LARGE \bf The equivariant Orlik-Solomon algebra}
\bigskip\\
{\bf Nicholas Proudfoot} \\
Department of Mathematics, University of California,
Berkeley, CA 94720
\bigskip
{\small
\begin{quote}
\noindent {\em Abstract.}
Given a real hyperplane arrangement $\A$, the complement $\MA$ of the complexification
of $\A$ admits an action of $\Zt$ by complex conjugation.
We define the equivariant Orlik-Solomon algebra of $\A$
to be the $\Zt$-equivariant cohomology ring
of $\MA$ with coefficients in $\Zt$.  We give a combinatorial presentation
of this ring, and  
interpret it as a deformation of the ordinary Orlik-Solomon algebra
into the Varchenko-Gel$'$fand ring of locally constant $\Zt$-valued
functions on the complement $\CA$ of $\A$ in $\R^n$.
We also show that the $\Zt$-equivariant homotopy type of $\MA$
is determined by the oriented matroid of $\A$.
As an application, we give two examples of pairs of arrangements
$\A$ and $\A'$ such that $\MA$ and $\mathcal{M}(\A')$ have the 
same nonequivariant homotopy
type, but are distinguished by the equivariant
Orlik-Solomon algebra.  
\end{quote}
}
\bigskip

\begin{section}{Introduction}\label{intro}
Let $\A = \{H_1,\ldots,H_n\}$ be an arrangement of $n$ hyperplanes in $\C^d$,
with $H_i = \omega_i^{-1}(0)$ for some affine linear map $\omega_i:\C^d\to\C$.
Let $\MA$ denote the complement of $\A$ in $\C^d$.
It is a fundamental problem in the study of hyperplane arrangements
to study the extent to which the topology of $\MA$ is determined
by the combinatorics of $\A$.

Let $\conea$ denote the central arrangement of hyperplanes in $\C^{d+1}$
given by first adding a ``hyperplane at infinity'' to $\A$ to produce an arrangement
of hyperplanes in $\C P^d$, and then taking its cone.
The {\em pointed matroid} of $\A$ is defined to be the matroid of dependence
relations among the hyperplanes of $\conea$, along with a specified basepoint
corresponding to the cone over the hyperplane at infinity \cite{F2}.
Geometrically, the pointed matroid encodes two types of data:
\begin{enumerate}
\item which subsets $S\subs\{1,\ldots,n\}$ have the property that
$\bigcap_{i\in S}H_i = \emptyset$, and
\item which subsets $S\subs\{1,\ldots,n\}$ have the property that
$\operatorname{codim}\bigcap_{i\in S}H_i > |S|$.
\end{enumerate}

\begin{definition}\label{osdef}
The {\em Orlik-Solomon algebra} $A(\A;R)$ is the cohomology ring
$H^*(\MA;R)$ of the complement of the complexified arrangement
with coefficients in the ring $R$.
\end{definition}

For each $i\leq n$, let $e_i = \omega_i^*[\R^+] \in A(\A;R)$
be the pullback of the generator $[\R^+]\in H^1(\C^*;R)$
under the map $\omega_i:\MA\to\C^* = \C\setminus\{0\}$.
The following theorem, due to Orlik and Solomon, states
that the elements $e_1,\ldots,e_n$ generate $A(\A;R)$, and
gives explicit relations in terms of the pointed matroid of $\A$.
We give here a simplified version by working only with the coefficient
ring $R=\Zt$, because this is the version 
that will extend well to the equivariant
setting.

\begin{theorem}\label{os}{\em\cite{OT}}
Consider the linear map 
$\partial = \sum_{i=1}^n\frac{\partial}{\partial e_i}$
from $\Zt[e_1,\ldots,e_n]$ to itself, lowering degree by 1.
The Orlik-Solomon algebra $A(\A;\Zt)$ is isomorphic to 
$\Zt[e_1,\ldots,e_n]\big/\mathcal{I}$,
where $\mathcal{I}$ is generated by the following three families
of relations:
\begin{eqnarray*}&1)&\hs\hs e_i^2\hs\hs\text{ for }\hs\hs i\in\{1,\ldots,n\}\\&&\\
&2)&\hs\hs
\prod_{i\in S}e_i
\hs\hs\text{ if }\hs\hs \bigcap_{i\in S}H_i = \emptyset\\&&\\
&3)&\hs\hs
\partial\prod_{i\in S}e_i
\hs\hs\text{ if }\hs\hs \bigcap_{i \in S}H_i
\text{ is nonempty with codimension greater than $|S|$.}
\end{eqnarray*}
\end{theorem}

Now suppose that our arrangement $\A$ is defined over the real numbers.
More precisely, suppose that $\omega_i$ restricts to a map $\omega_i:\R^d\to\R$
for all $i$.  
Let $$\Hip = \{p\mid \omega_i(p) > 0\}
\hspace{10pt}\text{     and     }\hspace{10pt}\Him = \{p\mid \omega_i(p) < 0\},$$
both {\em open} half-spaces in $\R^d$ with boundary $H_i$.
The {\em pointed oriented matroid} of $\A$ is defined to be the
oriented matroid with basepoint given by the dependence relations of $\conea$.
Like the pointed matroid, the pointed oriented matroid 
also encodes two types of geometrical data:
\begin{enumerate}
\item which subsets $S\subs\{1,\ldots,n\}$ have the property that
$\bigcap_{i\in S}H_i = \emptyset$, and
\item which pairs of subsets $S^+,S^-\subs\{1,\ldots,n\}$ have the property that
$\bigcap_{i\in S^+}H_i^+\cap\bigcap_{j\in S^-}H_j^- = \emptyset$.
\end{enumerate}
Note that this data is stronger than the data of the oriented matroid;
if $\cap_{i\in S} H_i$ is nonempty, then its codimension is greater than
$|S|$ if and only if there exists a decomposition
$S = \Sp\cup\Sm$ 
satisfying the second condition above.

In this paper we study the action of $\Zt=Gal(\C /\R)$ on $\MA$ by complex conjugation,
with fixed point set $\CA\subs\R^d$ equal to the complement of the real
loci of the hyperplanes.  This is an enhancement of the topological data of $\A$,
just as the pointed oriented matroid is an enhancement of the combinatorial data.
It is therefore natural to make the following definition.

\begin{definition}\label{eos}
The {\em equivariant Orlik-Solomon algebra} $\A_2(\A,\Zt)$ of a hyperplane
arrangement defined over $\R$ is the equivariant cohomology ring
$\hz(\MA;\Zt)$.
\end{definition}

In Section \ref{last} we give a presentation of the equivariant
Orlik-Solomon algebra in terms of the pointed oriented matroid of $\A$,
analogous to Theorem \ref{os}.\footnote{A special case of this presentation
first appeared in \cite[5.5]{HP}, using the geometry of hypertoric varieties.}
Moreover, we interpret $\A_2(\A,\Zt)$
as a deformation from the ordinary Orlik-Solomon algebra $A(\A;\Z)$
to the {\em Varchenko-Gel$\, '\!$fand ring} $\vgzt$, which is defined
to be the ring of locally constant functions from $\CA$ to $\Zt$.
We thus recover by independent means the presentation of $\vgzt$
given in \cite{VG}, and provide a topological explanation for the parallels that
Varchenko and Gel$'$fand observe between the the rings $A(\A;\Z)$
and $\vgzt$.
Note that, while the Orlik-Solomon algebra is super-commutative
and the Varchenko-Gel$'$fand ring is commutative, these two notions agree
in characteristic $2$.

A celebrated theorem of Salvetti \cite{Sa} says that if $\A$ is central and essential,
then $\MA$ is homotopy equivalent to a simplicial complex that
can be constructed from the oriented matroid\footnote{If $\A$ is central,
the oriented matroid and pointed oriented matroid encode the same data.} 
of $\A$ (see \cite{Sa},
\cite{Pa}, and \cite{GR}).  In Section \ref{salvetti}, we show that
this simplicial complex has a natural, combinatorially defined
action of $\Zt$, and that the homotopy equivalence is equivariant with
respect to this action.  Hence the oriented matroid of $\A$
in fact determines the equivariant homotopy type of $\MA$.
This observation provides an explanation for the recent discovery of Huisman
that the equivariant fundamental group of a line arrangement is determined
by its pointed oriented matroid \cite{Hu}.

In Example \ref{falk} we consider the famous first example 
of two real arrangements
with different pointed matroids, but with homotopy equivalent complements
\cite{F1}.  We show that these two arrangements are 
distinguished by the equivariant Orlik-Solomon algebra, hence the homotopy equivalence
cannot be made equivariant.
In Example \ref{vertical}, we consider two arrangements whose pointed oriented
matroids are related by a flip \cite{F1}.  This implies that their
complements are homotopy equivalent, and that their unoriented pointed matroids
are isomorphic, but once again their equivariant homotopy
types are distinguished by the equivariant Orlik-Solomon algebra.
We conclude with a conjecture regarding the relationship
between the combinatorial data and the equivariant 
topology of a real arrangement.
\newline

\noindent {\em Acknowledgments.}
The author is grateful to Graham Denham for pointing out the
similarity between the rings described in \cite{HP} and \cite{VG},
and to Michael Falk and David Speyer for their 
help in understanding many examples.
\end{section}

\begin{section}{Equivariant cohomology}\label{EC}
In this section we review some basic definitions and results from \cite{Bo}.
Let $X$ be a topological space equipped with an action of a group $G$.

\begin{definition}\label{xg}
Let $EG$ be a contractible space with a free $G$-action.
Then we put $$X_G := X \times_G EG = (X\times EG)/G$$
(well-defined up to homotopy equivalence), and
define the
$G$-equivariant cohomology of $X$
$$H^*_G(X) := H^*(X_G).$$
\end{definition}

The $G$-equivariant map from $X$ to a point induces a map
on cohomology in the other direction, hence $H^*_G(X)$ is a module
over $H^*_G(pt) \cong H^*(BG)$, where $BG=EG/G$ is the classifying space
for $G$.  Indeed, $H^*_G$ is a contravariant functor
from the category of $G$-spaces to the category of $H^*_G(pt)$-modules.

\begin{example}\label{rpi}
If $G=\Zt$, then we may take $EG = S^{\infty}$
and $BG = S^{\infty}/\Zt = \R P^{\infty}$.
Then $\hz(pt;\Zt) = H^*(\R P^{\infty};\Zt) \cong \Zt[x]$.
\end{example}

Suppose that $X$ is a finite-dimensional manifold, and
let $Y\subs X$ be a $G$-invariant submanifold.
We denote by $[Y] \in H^*_G(X)$ the cohomology class represented 
in Borel-Moore homology by the
finite-codimension submanifold $Y_G \subs X_G$.  
This will be our principal means
of understanding specific equivariant cohomology classes in this paper.
We will need two technical theorems about equivariant
cohomology,
both of which we state below.
Let $X$ be a $\Zt$-space, and let $F\subs X$ be the fixed point set.

\begin{theorem}\label{borel}{\em\cite[\S XII, 3.5]{Bo}}
Suppose that $F$ is nonempty,
the induced action of $\Zt$ on $H^*(X;\Zt)$
is trivial, and $H^*(X;\Zt)$ is generated in degree $1$.
Then the Leray-Serre spectral sequence for the fiber bundle
$X\hookto X_{\Zt}\to \R P^{\infty}$ collapses at the $E_2$ term.
\end{theorem}

\begin{corollary}\label{free}
Under the hypotheses of Theorem \ref{borel}, 
any additive basis from $H^*(X;\Zt)$ lifts to a $\Zt[x]$-basis
for $\hz(X;\Zt)$ (and any set of lifts will do).  In particular, $\hz(X;\Zt)$
is a free module over $\Zt[x]$.
\end{corollary}

\begin{theorem}\label{b2}{\em\cite[\S IV, 3.7(b)]{Bo}}
The restriction map $\hz(X;\Zt)\to\hz(F;\Zt)\cong H^*(F;\Zt)[x]$ 
is an isomorphism
in all degrees greater than the dimension of $X$.
\end{theorem}

Corollary \ref{free} says that we may interpret $\hz(X;\Zt)$ as a flat
family of rings over the $\Zt$ affine line.
The following corollary says that this family is a deformation
of $H^*(X;\Zt)$ into $H^*(F;\Zt)$.

\begin{corollary}\label{zero}
Under the hypotheses of Theorem \ref{borel}, 
$$H^*(X;\Zt)\cong \hz(X;\Zt)/\langle x\rangle$$
and 
$$H^*(F;\Zt)\cong \hz(X;\Zt)/\langle x-1\rangle.$$
\end{corollary}

\begin{proof}
The first statement follows immediately from Corollary \ref{free}.
For the second statement, consider the ring
$\hz(X;\Zt)[x^{-1}]$ obtained by formally inverting $x$.
Theorem \ref{b2} tells us that the restriction map
$$\hz(X;\Zt)[x^{-1}]\to\hz(F;\Zt)[x^{-1}]\cong H^*(F;\Zt)[x,x^{-1}]$$
is an isomorphism in high degree.  But this map commutes
with multiplication by $x$ and $x^{-1}$, so it must be an isomorphism
in every degree.
Setting $x$ equal to $1$, we obtain the desired result.
\end{proof}

The following example will be fundamental to our applications.

\begin{example}\label{cstar}
Let $X = \C^*$, with $\Zt$ acting by complex conjugation.
Since $X$ deformation-retracts equivariantly onto the compact
space $S^1$, Theorem \ref{borel} applies.
The image of $x$ under the standard map $\Zt[x]=\hz(pt,\Zt)\to\hz(X;\Zt)$
is the $\Zt$-equivariant Euler class of the topologically trivial
real line bundle with a nontrivial $\Zt$ action.  This bundle has
a $\Zt$-equivariant section, transverse to the zero section, vanishing
exactly on the real points of $X$, and is therefore represented
by the submanifold $\R^*\subs\C^*$.
Abusing notation, we will write $x=[\R^*]\in\hz(X;\Zt)$.
Let $y=[\R^+]\in\hz(X;\Zt)$.  
Then $x-y$ is represented by $\R^-$, therefore $y(x-y) = 0$.
Corollary \ref{free} says that $\hz(X;\Zt)$ is additively generated
by $x$ and $y$.
Since $\Zt[x,y]/\langle y(x-y)\rangle$ is already a free module
of rank $2$ over $\Zt[x]$, Corollary \ref{free} 
tells us that there can be no more relations.
\end{example}
\end{section}

\begin{section}{The equivariant Orlik-Solomon algebra}\label{last}
We now give a combinatorial presentation of the equivariant
Orlik-Solomon algebra.

\begin{theorem}\label{eq}
The ring
$A_2(\A;\Zt)$ is isomorphic to $\Zt[e_1,\ldots,e_n,x]\big/\mathcal{J}$,
where $\mathcal{J}$ is generated by the following three families
of relations:\,\footnote{Note that all of these relations are polynomial;
the $x^{-1}$ in the third relation cancels.}
\begin{eqnarray*}&1)&\hs\hs e_i(x-e_i)\hs\hs\text{ for }\hs\hs i\in\{1,\ldots,n\}\\&&\\
&2)&\hs\hs
\prod_{i\in\Sp}e_i\times\prod_{j\in\Sm}(x-e_j)
\hs\hs\text{ if }\hs\hs \bigcap_{i\in\Sp}\Hip\cap\bigcap_{j\in\Sm}H_j^- = \emptyset\\&&\\
&3)&\hs\hs
x^{-1}\left(\prod_{i\in\Sp}e_i\times\prod_{j\in\Sm}(x-e_j)
-\prod_{i\in\Sp}(x-e_i)\times\prod_{j\in\Sm}e_j\right)\\
&&\hspace{15pt}\text{ if }\hs\hs 
\bigcap_{i\in\Sp}\Hip\cap\bigcap_{j\in\Sm}H_j^-=\emptyset
\hs\hs\text{ and }\hs\hs\bigcap_{i\in S}H_i \neq\emptyset, \hs\hs\text{ where }
\hs\hs S=\Sp\cup\Sm.
\end{eqnarray*}
\end{theorem}

\begin{proof}
Let $y = [\R^+] \in \hz(\C^*;\Zt),$ and let $$e_i = \omega_i^*(y) \in A_2(\A;\Zt),$$
represented by the submanifold $$Y_i^+ = \omega_i^{-1}(\R^+).$$
Let $x\in A_2(\A;\Zt)$ be the image of the generator of $\hz(pt;\Zt)$;
by functoriality, we have $x = \omega_i^*(x)$ for all $i$.
Recall from Example \ref{cstar} that $[\R^-] = x-y \in \hz(\C^*;\Zt)$,
hence $$x-e_i = \omega^*(x-y)\in A_2(\A;\Zt)$$ is represented by the submanifold
$$Y_i^- = \omega^{-1}(\R^-).$$
Theorem \ref{os} tells us that $e_1,\ldots,e_n$ are lifts of ring
generators for the ordinary Orlik-Solomon algebra $A(\A;\Zt)$.
Since the manifolds
$Y_i^+$ are stable under the action of $\Zt$, the induced action of $\Zt$
on $A(\A;\Zt)$ is trivial.
The space $\MA$ has a compact $\Zt$-equivariant deformation retract,
therefore Corollary \ref{borel} tells us that $ A_2(\A;\Zt)$ 
is generated as a ring
by the classes $e_i$ and $x$.  We must now check that each of the three families
of generators of $\mathcal{J}$ do indeed vanish in $ A_2(\A;\Zt)$, and that
they generate all of the relations.

The first family of relations follows from the fact that $Y_i^+\cap Y_i^-=\emptyset$
for all $i\in\{1,\ldots,n\}$.
For the second family, we must show that if 
$$\bigcap_{i\in\Sp}\Hip\cap\bigcap_{j\in\Sm}H_j^- = \emptyset,$$
then $$\bigcap_{i\in\Sp}Y_i^+\cap\bigcap_{j\in\Sm}Y_j^- = \emptyset.$$
Suppose that $$p\in\bigcap_{i\in\Sp}Y_i^+\cap\bigcap_{j\in\Sm}Y_j^-,$$
in other words $\omega_i(p)\in\R^+$ for all $i\in S^+$ and
$\omega_j(p)\in\R^-$ for all $j\in S^-$.
Then the real part $$\Re(p)\in\bigcap_{i\in\Sp}\Hip\cap\bigcap_{j\in\Sm}H_j^-,$$
hence the intersection is not empty.

For the third family, note that 
since $A_2(\A;\Zt)$ is free over $\Zt[x]$, it is sufficient
to show that 
\begin{equation}\label{two}
\prod_{i\in\Sp}e_i\times\prod_{j\in\Sm}(x-e_j)
-\prod_{i\in\Sp}(x-e_i)\times\prod_{j\in\Sm}e_j = 0.
\end{equation}
By hypothesis, $$\bigcap_{i\in\Sp}\Hip\cap\bigcap_{j\in\Sm}H_j^- = \emptyset
\hs\hs\text{ and }\hs\hs\bigcap_{i\in S}H_i\neq\emptyset,$$ 
where $S = \Sp\cup\Sm$.
Choose a point $p\in\bigcap_{i\in S}H_i$.  The involution
of $\R^d$ given by reflection through $p$ takes $\Hip$ to $\Him$
and vice versa for all $i\in S$, hence we also have
$$\bigcap_{i\in\Sp}\Him\cap\bigcap_{j\in\Sm}H_j^+ = \emptyset.$$
Thus both of the terms in Equation \eqref{two} are in fact
contained in the second family of relations, and are therefore
equal to zero in the ring $A_2(\A;\Zt)$.

Now we must show that we have found all of the relations.
Let $$\psi:\Zt[e_1,\ldots,e_n,x]\to\Zt[e_1,\ldots,e_n]$$
be the map given by sending $x$
to zero, and note that $\psi(\mathcal{J}) = \mathcal{I}$.
Now suppose that $\a\in\Zt[e_1,\ldots,e_n,x]$ is a relation in $ A_2(\A;\Zt)$ that is
{\em not} in the ideal $\mathcal{J}$, and choose $\a$ of minimal degree.  
By Corollary \ref{zero} we must have
$\psi(\a)\in\mathcal{I}$, hence there exists $\beta\in\mathcal{J}$
with $\psi(\a-\beta)=0$.  This implies that $\a-\beta = x\gamma$ for some
$\gamma\in\Zt[e_1,\ldots,e_n,x]$. 
Since $\a$ and $\beta$ are both relations in $A_2(\A;\Zt)$, 
and $A_2(\A;\Zt)$ is free over $\Zt[x]$,
$\gamma$ must also be a relation.  Since $\beta$ is in $\mathcal{J}$ and $\a$ is not,
$\gamma$ cannot be in $\mathcal{J}$.  Since $\deg\gamma = \deg\a - 1$, 
we have reached a contradiction. 
\end{proof}

By Corollary \ref{zero}, $A_2(\A;\Zt)$ is a flat family of rings parameterized
by the affine line $\operatorname{Spec}\Zt[x]$, specializing at $x=0$
to $H^*(\MA;\Zt)=A(\A;\Zt)$, and at $x=1$ to $H^*(\CA;\Zt)=\vgzt$.
In particular, this provides a topological explanation for the fact that
the dimension of the Orlik-Solomon algebra is equal to the number
of connected components of $\CA$.  
By setting $x=1$ in Theorem \ref{eq} we obtain a nontrivial presentation
of $\vgzt$, first given in the central case (over the integers) in \cite{VG}.
Varchenko and Gel$'$fand 
interpret $e_i\in\vgzt$ as the $i^{\text{th}}$ {\em Heaviside function}
$\CA\to\R$, restricting to $1$ on $\CA\cap H_i^+$ and $0$ on $\CA\cap H_i^-$.
These functions are easily seen to generate the ring $\vgzt$, and the three
families of relations are clear, but the proof that there are no other
relations is nontrivial.
Varchenko and Gel$'$fand observe that this presentation defines a filtration on
$\vgzt$ with $A(\A;\Zt)$ as its associate graded.  This is also a consequence
of Corollaries \ref{free} and \ref{zero}; this phenomenon is explored in greater
detail in \cite{Ca}.

\begin{remark}
Our presentations of $\vgzt$ and $A_2(\A;\Zt)$ depend on the
coorientations of the hyperplanes, while the isomorphism classes of
the rings themselves do not.
Reversing the orientation of the hyperplane $H_i$ corresponds to 
changing every appearance of $e_i$ to $x-e_i$ in the generators of
$\mathcal{J}$.
\end{remark}
\end{section}

\begin{section}{The Salvetti complex}\label{salvetti}
Let $\A$ be an essential central arrangement in $\R^d$.
Salvetti \cite{Sa} has constructed a simplicial complex from a poset $\Sala$,
depending only on the oriented matroid of $\A$, which is homotopy
equivalent to the complement $\MA$ of the complexification of $\A$.
In this section we define a combinatorial action of $\Zt$ on $\Sala$,
and show that the homotopy equivalence is equivariant.

The hyperplanes of $\A$ subdivide $\R^d$ into faces, open in 
their supports, which form a poset
$\F$ ordered by reverse inclusion.  
The minimal elements of $\F$ are the connected components of $\CA$, and
$\{0\}$ is the unique maximal element.
The {\em Salvetti poset} $\Sala$ is a poset consisting of elements of the form
$$\{(F,C)\mid C\text{ minimal and }C\leq F\}.$$
The partial order is determined by putting $(F',C')\leq (F,C)$ if and only if 
$F'\leq F$ and $C' = F'C$, where the latter equality means 
that $C$ and $C'$ lie on the same side of every hyperplane containing $F'$.
The {\em Salvetti complex} $|\Sala|$ is defined to be the order complex of this poset.

The poset $\Sala$ admits an action of $\Zt$
given by setting $(F,C)^* = (F,\tilde{C})$, where $\tilde{C}$ is obtained from $C$
by reflecting it over all of the hyperplanes that contain $F$.
In \cite{GR}, $\Sala$ is defined as a subset
of the set of all functions from the ground set of the oriented matroid to
the set $\{\pm 1, \pm i\}$.  In this language, our $\Zt$-action is simply
complex conjugation, and is easily seen to be
an invariant of the oriented matroid.
This action induces an action of $\Zt$ on the Salvetti complex $|\Sala|$.

\begin{theorem}\label{eqsal}
The complex $|\Sala|$ is equivariantly homotopy equivalent to $\MA$.
In particular, the equivariant homotopy type of $\MA$ is determined
by the oriented matroid associated to $\A$.
\end{theorem}

\begin{proof}
For every $F\in\F$, choose a point $x(F)\in F\subs\R^d$.
Each element of $\Sala$ determines a vertex in the complex $|\Sala|$.
For all $(F,C)\in\Sala$, let 
$$V(F,C) = 
\begin{cases}
\left\{{\displaystyle\sum_{C'\leq F}}\la_{C'}x(C')\hs\hs\bigg |\hs\hs\la_{C'}>0\right\} &\text{if }F\neq\{0\}\\
\R^d &\text{if }F=\{0\},
\end{cases}$$
and let
$$W(F,C) = \{x\in\R^d\mid x\text{ and }C\text{ lie on the same side
of every hyperplane containing }F\}.$$
Paris \cite{Pa} shows that 
$$\mathcal{U} = \Big\{V(F,C) + iW(F,C)\hs\hs\Big |\hs\hs (F,C)\in\Sala\Big\}$$
is an open cover of $\MA$ with nerve $|\Sala|$, and that any nonempty
intersection of open sets from $\mathcal{U}$ is contractible,
hence concluding that $\MA$ is homotopy equivalent to $|\Sala|$.
To extend this proof to the equivariant context, we need only
observe that $W(F,\tilde{C}) = W(F,C)$,
and $V(F,\tilde{C}) = -V(F,C)$.  Both of these equalities are clear
from the definitions.
\end{proof}

\vspace{-\baselineskip}
\begin{remark}
The Salvetti complex may be defined for an arbitrary oriented matroid,
which may not be realizable by a hyperplane arrangement (see for example
\cite{BLSWZ}.  We can then define the equivariant Orlik-Solomon algebra
of an arbitrary oriented matroid to be the $\Zt$-equivariant cohomology
ring of its Salvetti complex.  Theorem \ref{eqsal} implies that this definition
agrees with our original one if the oriented matroid is realizable.
\end{remark}
\end{section}

\begin{section}{Examples}\label{examples}
In this section we discuss two examples in which
the equivariant Orlik-Solomon algebra
distinguishes two arrangements with (nonequivariantly)
homotopy equivalent complements.  
In both, we work with affine arrangements to keep dimensions
as low as possible.  The analogous central examples
can be understood via the following proposition.

\begin{proposition}\label{cone}
There is a $\Zt$-equivariant diffeomorphism
$\mathcal{M}(\mathcal{CA}) \cong\MA\times\C^*$, and
$$A_2(\mathcal{CA};\Zt) \cong A_2(\A;\Zt)\otimes_{\Zt[x]}\Zt[x,y]/y(x-y).$$
\end{proposition}

\begin{proof}
The standard diffeomorphism 
$\mathcal{M}(\mathcal{CA}) \cong\MA\times\C^*$, found for example in \cite{OT},
is $\Zt$-equivariant.  The second half of the proposition 
is simply the statement
of the equivariant K\"unneth theorem \cite[7.4]{Se}, 
combined with Example \ref{cstar}.
\end{proof}

\vspace{-\baselineskip}
\begin{example}\label{falk}
The example of Figure \ref{cool} was introduced by Falk \cite[3.1]{F1}.
\begin{figure}[h]
\begin{center}
\psfrag{A}{$\A$}
\psfrag{A'}{$\A'$}
\includegraphics[height=1.5in]{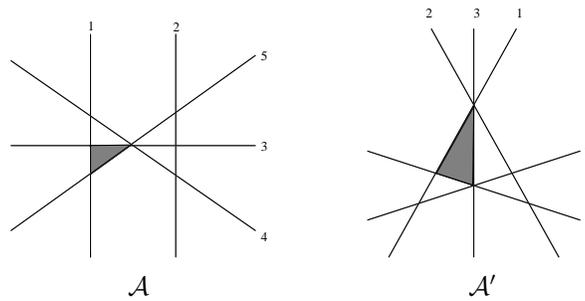}
\caption{Two arrangements whose complements are homotopy equivalent,
but not equivariantly.}
\label{cool}
\end{center}
\end{figure}
The arrangements $\A$ and $\A'$ have nonisomorphic 
pointed matroids,
but their complements are homotopy equivalent.  In particular,
they cannot be distinguished by their Orlik-Solomon algebras.
We show that their equivariant Orlik-Solomon algebras are nonisomorphic,
therefore the homotopy equivalence between their complements cannot
be $\Zt$-equivariant.
Choose coorientations so that the intersections $\cap_{i\leq 5}\Him$
are equal to the shaded regions.  Then
$$A_2(\A;\Zt)=\Zt[e_1,\ldots,e_5,x]/\mathcal{J}\hspace{10pt}\text{ and }
\hspace{10pt}A_2(\A';\Zt)=\Zt[e_1,\ldots,e_5,x]/\mathcal{J'},$$
where $$\mathcal{J} = \left<
\begin{array}{c}
e_1(x-e_1),\ldots,e_5(x-e_5), 
e_1e_2, e_1(x-e_3)e_4,
e_1e_3e_5, e_1e_4e_5, 
e_2e_3(x-e_4),\\
e_2(x-e_4)(x-e_5),
e_2(x-e_3)(x-e_5),
e_3e_4+e_3e_5+e_4e_5+e_4x
\end{array}
\right>$$
and 
$$\mathcal{J'} = \left<
\begin{array}{c}
e_1(x-e_1),\ldots,e_5(x-e_5), 
e_1e_2e_4, e_1e_2e_5, e_1e_3e_4, e_1e_3e_5,
e_1e_4(x-e_5), e_2(x-e_3)e_4,\\
e_2(x-e_3)e_5, e_2(x-e_4)e_5, 
e_1e_2+e_1e_3+e_2e_3+e_2x,
e_3e_4+e_3e_5+e_4e_5+e_4x
\end{array}
\right>.$$
Using Macaulay 2 \cite{M2}, we can compute the number of elements
in each degree that annihilate a given linear form.  By comparing the
lists that we obtain for the two rings, we conclude that the rings are not isomorphic.

These two arrangements are generic rank $2$ truncations of a pair
of rank $3$ arrangements $\A_3$ and $\A_3'$
which have diffeomorphic complements by a
general construction relating parallel connections to direct sums
(see \cite[Thm 2]{EF} and \cite[3.8]{F2}).  
The first arrangement $\A_3$ is given by the equation
$(x+1)(x-1)y(y+z)(y-z)=0$, with $\A$ obtained from $\A_3$ by setting $z=x$.
The second arrangement $\A_3'$ is given by the equation
$(2x+y-z)(2x-y+z)x(x-y)(x+y)=0$, with $\A'$ obtained 
from $\A_3'$ by setting $z=1$.
The diffeomorphism between $\mathcal{M}(\A_3)$ and $\mathcal{M}(\A_3')$
given in \cite{EF} is easily seen to be $\Zt$-equivariant, 
as it is essentially
derived from repeated applications of the diffeomorphism
of Proposition \ref{cone}.
Furthermore, it is not hard to produce an explicit isomorphism
between $A_2(\A_3;\Zt)$ and $A_2(\A_3';\Zt)$.  This shows that
a theorem of Pendergrass \cite[3.11]{F2}, which states that
truncation of matroids preserves isomorphisms of Orlik-Solomon
algebras, does not extend to the equivariant setting.
\end{example}

\begin{example}\label{vertical}
Consider the two line arrangements shown in Figure \ref{save}.\footnote{This
example appeared first in \cite{HP}.}
\begin{figure}[h]
\begin{center}
\psfrag{A}{$\A$}
\psfrag{A'}{$\A'$}
\includegraphics{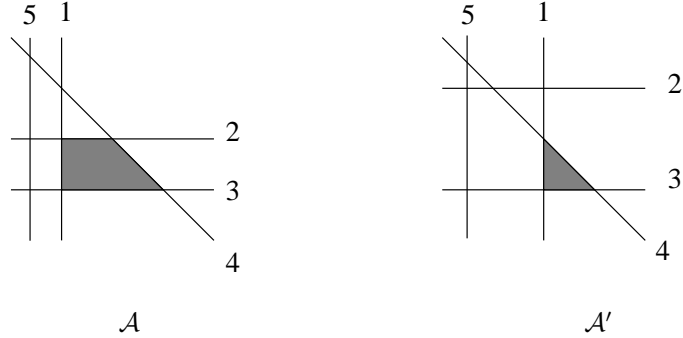}
\caption{Two arrangements related by a flip
with nonisomorphic Orlik-Solomon algebras.}
\label{save}
\end{center}
\end{figure}

The pointed oriented matroids corresponding to these two arrangements
are related by a flip.  Geometrically, this means that
$\A'$ can be obtained from $\A$ 
by translating one of the hyperplanes from one side
of a vertex to another.  (For a precise definition of flips, 
see \cite[\S 7.3]{BLSWZ}.)
Falk \cite{F1} has shown that any two real line arrangements related by a flip
have homotopy equivalent complements; in this example we show 
that Falk's theorem does not extend to the equivariant setting.
We have
$$A_2(\A;\Zt)\cong
\Zt [\vec{e},x]\bigg/
\< \begin{array}{c}
e_1(x-e_1), e_2(x-e_2), e_3(x-e_3), e_4(x-e_4),\\ e_5(x-e_5), 
e_2e_3, (x-e_1)e_5, e_1(x-e_2)e_4,\\ e_1e_3e_4, (x-e_2)e_4e_5, e_3e_4e_5 
\end{array}\>$$
and
$$A_2(\A';\Zt)\cong
\Zt [\vec{e},x]\bigg/
\< \begin{array}{c}
e_1(x-e_1), e_2(x-e_2), e_3(x-e_3), e_4(x-e_4),\\ e_5(x-e_5),
e_2e_3, (x-e_1)e_5, (x-e_1)e_2(x-e_4),\\ e_1e_3e_4, (x-e_2)e_4e_5, e_3e_4e_5 
\end{array}\>.$$
Using Macaulay 2 \cite{M2}, we find that the annihilator of the element
$e_2\in A_2(\A;\Zt)$ is generated by two linear elements
(namely $e_3$ and $x-e_2$) and nothing else, while none of the
(finitely many) elements
of $A_2(\A';\Zt)$ has this property.
Hence the two rings are not isomorphic, and $\MA$ is not equivariantly
homotopy equivalent to $\mathcal{M}(\A')$.
From this example we conclude that the equivariant Orlik-Solomon algebra
of an arrangement is {\em not} determined by the pointed 
{\em un}oriented matroid.
\end{example}

Both of the arrangements that we have discussed
have connected pointed matroids.
Eschenbrenner and Falk \cite{EF} conjecture that if $\A$ is
a complex central arrangement with connected matroid, then the matroid
of $\A$ is determined by the homotopy type of $\MA$.
The analogous conjecture for real arrangements may be stated as follows.

\begin{conjecture}
If $\A$ is
a real central arrangement with connected matroid, then the oriented matroid
of $\A$ is determined by the equivariant homotopy type of $\MA$.
\end{conjecture}
\end{section}

\end{spacing}
\end{document}